\documentclass[11pt]{article}

\def\BC#1#2{\ensuremath{\left({#1\atop#2}\right)}}  
\def\eqref#1{\mbox{(\ref{eq:#1})}}    

\begin{document}

\title{Erlang's loss function \\ and \\ a property of the exponential function}
\author{Hans J.H. Tuenter \\ \\
        Mathematical Finance Program,  \\
        University of Toronto, \\
        720 Spadina Avenue, Suite 219,  \\
        Toronto, Ontario,           \\
        Canada M5S 2T9. \\[0.25cm]
        Email: hans.tuenter@utoronto.ca}
\date{\small \begin{tabular}{ll}
               First Version: & August 2005 \\
               This Version:  & November 2006 
            \end{tabular}      
}

\maketitle
\begin{abstract}
 We show that a concavity property of the exponential function is a direct consequence of the convexity of
the continued Erlang loss function.
%
\vspace{-.25cm}
\paragraph*{Key~Words.} Exponential function, Erlang's loss function, Convexity.
\end{abstract}


\section{Introduction}
For positive $\lambda$, let $s_n(\lambda)$ denote the $n$th partial sum of the Taylor series for the exponential function: 
\begin{equation}
  s_n(\lambda)=\sum_{j=0}^n\frac{\lambda^j}{j!}.
 \label{eq:1.1}
\end{equation}
Berenstein et al.~\cite{BerensteinVK:1986} have shown that,
for nonnegative integers~$n$, 
this sum has the following property:
\begin{equation}
  \frac{s_{n}(\lambda)}{s_{n+1}(\lambda)} + 
  \frac{s_{n+2}(\lambda)}{s_{n+3}(\lambda)}  <
  2\frac{s_{n+1}(\lambda)}{s_{n+2}(\lambda)},
 \label{eq:1.2}
\end{equation} 
and, as a consequence, that
\begin{equation}
  \frac{\BC{n-m+l}{l}}{\BC{n}{l}}s_{n-l}(\lambda)s_{n-m+l}(\lambda) < 
  s_n(\lambda)s_{n-m}(\lambda)<s_{n-l}(\lambda)s_{n-m+l}(\lambda),
 \label{eq:1.3}
\end{equation}
holds for positive integers $n\ge m>l$.
These inequalities are then used to establish bounds for the mean length of the queue 
in the $E_k|\mathit{GI}|1|\infty$ system,
and to prove a monotonicity and convexity property of the mean number of customers in the first stage of a 
closed, two-stage $\mathit{GI}|M|1|N$ queueing system.
In a subsequent paper, Berenstein and Vainshtein~\cite{BerensteinVainshtein:1997} use~\eqref{1.2} and~\eqref{1.3} 
as a springboard to derive more general results.

We show that the above inequalities can be established in an elegant and succinct manner by using the 
result by Jagers and van Doorn on the convexity of the analytic
continuation of Erlang's loss function~\cite{JagersDoorn:1986}. 

\section{Erlang's loss function}
For positive $\lambda$ and nonnegative integers~$n$, Erlang's loss function is defined as
\begin{equation}
  B(n,\lambda)=\frac{\lambda^n}{n!}/\sum_{j=0}^n\frac{\lambda^j}{j!},
  \label{eq:2.1}
\end{equation}
and represents the probability that calls arriving in a Poisson stream
of rate $\lambda$ to a link consisting of $n$ circuits and requesting
unit service time would find all circuits occupied and therefore be lost.
The analytic continuation of~\eqref{2.1} is given by
\begin{equation}
  B(x,\lambda)=\left\{\lambda\int_0^\infty e^{-\lambda
      t}(1+t)^x\,dt\right\}^{-1}.
  \label{eq:2.2}
\end{equation}
Integration by parts shows that
\begin{equation}
 B(x+1,\lambda)^{-1}=1+\frac{x+1}{\lambda}B(x,\lambda)^{-1}.
 \label{eq:2.3}
\end{equation}
Using this recursion, with initial value $B(0,\lambda)=1$, 
it is easily verified that~\eqref{2.2} also coincides with~\eqref{2.1}
for positive, integral values of~$x$.
Jagers and van Doorn~\cite{JagersDoorn:1986} have shown that the analytic continuation
$B(x,\lambda)$ is a strictly convex function 
of $x\ge0$, for every $\lambda>0$. 

\section{Derivation}
The ratio $\varphi(n,\lambda)=s_{n}(\lambda)/s_{n+1}(\lambda)$ can be expressed in terms of the Erlang loss function:
\[ \varphi(n,\lambda)
   =\frac{s_{n}(\lambda)}{s_{n+1}(\lambda)}=\frac{n+1}{\lambda}\frac{B(n+1,\lambda)}{B(n,\lambda)}
   =1-B(n+1,\lambda),
\]
where we have used~\eqref{2.3} to substitute for $B(n,\lambda)$.
Since $B(x,\lambda)$ is strictly convex, the analytic continuation
$\varphi(x,\lambda)$ is strictly concave, and~\eqref{1.2} 
is an easy consequence.
Moreover, this shows that, for nonnegative integers~$n$, we have
\begin{equation}
  \frac{m-r}{m}\varphi(n,\lambda)+\frac{r}{m}\varphi(n+m,\lambda)<\varphi(n+r,\lambda),
 \label{eq:3:2}
\end{equation}
with integers $1\le r<m$.
To prove the left-hand side of~\eqref{1.3}, note that $B(x,\lambda)$ is strictly convex for $x\ge0$, 
so that $\varphi(x,\lambda)$ is strictly concave for $x\ge-1$.
Since $B(0,\lambda)=1$, we can define $\varphi(-1,\lambda)=0$, and extend~\eqref{3:2}
to include the case~$n=-1$. This gives
\begin{equation}
  \frac{\varphi(m,\lambda)}{m+1}<\frac{\varphi(r,\lambda)}{r+1},
  \label{eq:3.3}
\end{equation}
with integers $0\le r<m$. 
In particular, the choice of $m=n=r+1$  
immediately gives the inequality $s_{n-1}(\lambda)s_{n+1}(\lambda)>\frac{n}{n+1}s_n^2(\lambda)$,
valid for $n\ge1$. More generally, for $1\le l<m\le n$, we have 
\[ \frac{(n-l)!}{n!}\frac{s_{n-l}(\lambda)}{s_{n}(\lambda)} 
   = \prod_{j=n-l}^{n-1}\frac{\varphi(j)}{j+1} < 
     \prod_{j=n-l}^{n-1}\frac{\varphi(l-m+j)}{l-m+j+1} 
   = \frac{(n-m)!}{(n-m+l)!}\frac{s_{n-m}(\lambda)}{s_{n-m+l}(\lambda)},
\]
where the central inequality follows from a repeated, pairwise application of~\eqref{3.3}. 
Cross multiplying the outer terms, and collecting the factorials in binomial coefficients,
establishes the left-hand side of~\eqref{1.3}.
The right-hand side of~\eqref{1.3} is a direct consequence of the fact that
$B(x,\lambda)$ is a strictly decreasing function of~$x$, as is easily seen from its definition~\eqref{2.2},
and implies that $\varphi(n,\lambda)$ is a strictly increasing function of~$n$.
This then gives
\[ \frac{s_{n-m}(\lambda)}{s_{l-m+l}(\lambda)}=\prod_{j=n-l}^{n-1}\varphi(l-m+j,\lambda)
   <\prod_{j=n-l}^{n-1}\varphi(j,\lambda)=\frac{s_{n-l}(\lambda)}{s_n(\lambda)},
\]
for $1\le l<m\le n$, and establishes the right-hand side of~\eqref{1.3}.

\section{Comment}
As we have shown, properties~\eqref{1.2} and~\eqref{1.3} of the sum $s_n(\lambda)$ 
are a direct consequence of the convexity of the continued Erlang loss function.
These properties can, of course, be proven by purely algebraic manipulations, as in~\cite{BerensteinVK:1986}.
However, it is gratifying to solve a number theoretical problem that originated in queueing theory
with a result from queueing theory.


\footnotesize

\begin{thebibliography}{1}

\bibitem{BerensteinVK:1986}
A.D. Berenshtein, A.D. Vainshtein, and A.Ya. Kreinin.
\newblock A convexity property of the {P}oisson distribution and its
  application in queueing theory.
\newblock {\em Journal of Mathematical Sciences}, 47(1):2288--2292, October
  1989.
\newblock Translated from Problemy Usto\v{\i}chivosti Stokhasticheskikh
  Modele{\v\i}, Trudy Seminara (Stability Problems of Stochastic Models,
  Seminar Proceedings), pp.~{17--22}, 1986.

\bibitem{BerensteinVainshtein:1997}
Arkady Berenstein and Alek Vainshtein.
\newblock Concavity of weighted arithmetic means with applications.
\newblock {\em Archiv der Mathematik}, 69(2):120--126, August 1997.

\bibitem{JagersDoorn:1986}
A.A. Jagers and Erik~A. van Doorn.
\newblock On the continued {E}rlang loss function.
\newblock {\em Operations Research Letters}, 5(1):43--46, June 1986.
\end{thebibliography}
\end{document}